\documentclass[12pt]
{amsart}
\usepackage{amsfonts,amsmath,amssymb,amsthm,graphicx,afterpage,url}
\usepackage{hyperref}
\usepackage{color}

\input{xy}
\xyoption{all}

\def\R{{\mathbb R}} \def\Z{{\mathbb Z}}

\long\def\comment#1\endcomment{}

\newcommand{\jonly}[1]{}

    \theoremstyle{theorem}
         \newtheorem{theorem}{Theorem}

    \theoremstyle{definition}
         \newtheorem{remark}[theorem]{Remark}

\begin{document}

\title{To S. Parsa's theorem on embeddability of joins}

\author{A. Skopenkov}

\thanks{Independent University of Moscow, Moscow Institute of Physics and Technology, Russian Federation.
\texttt{https://users.mccme.ru/skopenko}.
}

\date{}

\maketitle



The purpose of this short note is to justify (in Remark \ref{r:just}) that \cite{Pa21} is not a reliable reference for its main result (the part most interesting to mathematicians not specialized in the area is restated as Theorem \ref{t:parsa} below).
We also present a reliable reference \cite[Theorem A and Remark 1.1]{Me22}.
This note is an extended version of the Zentralblatt review of \cite{Pa21}.

\smallskip
{\bf Disclaimer.} (a) For a discussion of what is a reliable reference see e.g. \cite[p. 2]{Sk21d}.
No other meaning of `not a reliable reference' is meant here.
In particular, I do not mean that the proofs of \cite{Pa21} cannot be corrected; I only mean that their correction is the author's job not the reader's job, see \cite[Remark 2.3]{Sk21d}.\footnote{I am sorry S. Parsa did not accept my proposal to revise \cite{Pa21} using Remark \ref{r:just}, and to update arXiv version, so that I could praise the latter instead of criticizing the journal version.
(Instead, he sent me some ideas of how \cite{Pa21} can be revised, cf. \cite[Remark 2.3]{Sk21d}.)
Open publication of criticism adds an example of different reliability standards in current mathematical research \cite{Sk21d}, and could stimulate appearance of a reliable revision of \cite{Pa21}.}

(b) No priority question is raised here.
Both the paper \cite{Me22} and the current paper attribute Theorem \ref{t:parsa} to S. Parsa.

\smallskip
We abbreviate `$k$-dimensional finite simplicial complex' to `$k$-complex'.

Clearly, if a complex $K$ embeds into $\R^m$, and a complex $L$ embeds into $\R^n$, then the join $K*L$  embeds into $\R^{m+n+1}$.
Under some conditions `non-embeddability of $k$-complexes into $\R^{2k}$ is stable under joins', i.e.,
if $K, L$ are complexes of dimensions $k,l$ such that $K$ does not embed into $\R^{2k}$, and $L$ does not embed into $\R^{2l}$, then $K*L$ does not embed into $\R^{2(k+l+1)}$, see \cite{BKK, MS06, Pa20, PS20, Me22}.

\begin{theorem}[Parsa]\label{t:parsa}
For any $k,l\ge2$ there exist complexes $K, L$ of dimensions $k,l$ such that $K$ does not embed into $\R^{2k}$,
and $L$ does not embed into $\R^{2l}$, but $K*L$  embeds into $\R^{2(k+l+1)}$.
\end{theorem}

This result is based on delicate difference between the integer-valued and the modulo 2 van Kampen obstructions.
A complicated (and unreliable, see Remark \ref{r:just}) proof is given in \cite{Pa21}.
A simpler proof given in \cite{Me22} uses `cobordism interpretation' of the van Kampen obstructions.
It would be interesting to further simplify the proof of \cite{Me22} by using the van Kampen obstructions themselves, not their `cobordism interpretation'.

\begin{remark}[on \cite{Pa21}]\label{r:just}
Note that p. 7175-1176 of the published version of \cite{Pa21} correspond pp. 23-24 of the arXiv v4 version.

(1) Proof of Theorem 1 in p. 7175 is incomplete.
Indeed, the join $M*N$ of Theorem 1 is not at all considered in that proof.

(1a) Less important remarks: in the last sentence of the proof of Theorem 1 in p. 7175
it is not explained and it is not clear,

$\bullet$ what is `the claim'.

$\bullet$ how the last sentence of the first paragraph follows from the sentence on the `sequence of moduli' of Lemma 10 (this last sentence follows from the `in particular' part of Lemma 10).

(2) Theorem 1 also follows from the case
{\it `if both the van Kampen obstructions of $M$ and $N$ vanish modulo 2'} of Theorem 7(1) in p. 7175.
However, this case is not proved in the proof of Theorem 7 in p. 7175-7176.
Indeed, the subcase
{\it `if both the van Kampen obstructions of $M$ and $N$ vanish modulo 2 and $d_M+d_N\ge2$'}
 is not at all considered there
(so
in p. 7176 the phrase {\it `Since the van Kampen obstruction of $M$ vanishes'} is incorrect, and the phrase {\it `this proves statement (1) in the case $d_M+d_N\ge2$'} is unjustified).

(3) Less important remark.
The main definition of $I_p(K)$ before Theorem 2 is meaningless to a reader, because it uses the notion
{\it `the Smith index of the $\Z_p$-complex $K$ computed with coefficients $\Z_p$'},
and for this notion neither definition nor reference to the definition is presented.



(4) The paper \cite{Pa21} is quite technical and (as above shows) poorly written.
So the paper is not accessible to a large audience, as opposed to
\cite[p. 7151, the paragraph before `Method of proof']{Pa21}.
In that paragraph the complicated method of proof of Theorem \ref{t:parsa} in \cite{Pa21}
is justified by the hope to find applications in other problems of combinatorics and discrete geometry.
However, no results (or even speculations) justifying this hope are presented.
Currently there are only two results \cite{Pa20, Pa21} obtained by this method of proof;
they have much simpler proofs \cite{PS20, Me22}; this information (necessary for the reader to decide whether to share the author's hope) is omitted in the paragraph expressing this hope.
(The paper \cite{PS20} was known to the author before publication of \cite{Pa21};
the paper \cite{Me22} could be unknown, but could be quoted in arXiv update of \cite{Pa21}.)
\end{remark}

\end{document}